\input amstex
\documentstyle{amsppt}
\magnification=\magstep1
\baselineskip=18pt
\hsize=6truein
\vsize=8truein

\topmatter
\title On the uniqueness of AdS space-time in higher dimensions \endtitle
\author Jie Qing \\
Department of Mathematics\\
UCSC
\endauthor

\leftheadtext{Asymptotically AdS}
\rightheadtext{Jie Qing}

\address Jie Qing, Dept. of Math., UC,
Santa Cruz, Santa Cruz, CA 95064.
\endaddress

\email qing{\@}math.ucsc.edu \endemail

\abstract In this paper, based on an intrinsic definition of
asymptotically AdS space-times, we show that the standard anti-de
Sitter space-time is the unique strictly stationary asymptotically
AdS solution to the vacuum Einstein equations with negative
cosmological constant in dimension less than $7$. Instead of using
the positive energy theorem for asymptotically hyperbolic spaces
our approach appeals to the classic positive mass theorem for
asymptotically flat spaces.
\endabstract

\endtopmatter

\document

\head 1. Introduction \endhead

Recently, there has been some interest in the study of space-times
that satisfy Einstein equations with negative cosmological
constant in association with the so-called AdS/CFT correspondence.
With the presence of a negative cosmological constant, the anti-de
Sitter space-time replaces the Minkowski space-time as the ground
state of the theory. Bocher, Gibbons and Horowitz showed that in
$3+1$ dimensions, the only strictly stationary asymptotically AdS
space-time that satisfies the vacuum Einstein equations with
negative cosmological constant is the anti-de Sitter space-time in
[BGH] (see also [CS]). Another class of globally static
asymptotically locally AdS space-times, the AdS solitons, are also
important in the theory. In [GSW], Galloway, Surya and Woolgar
proved a uniqueness theorem of the AdS solitons. Later, in [ACD],
Anderson, Chrusciel and Delay improved the uniqueness theorem of
AdS solitons.

In [Wa1], the uniqueness result of [BGH] was generalized to higher
dimensions when the space-time is static and the static slice is
of spin structure. Proofs in [BGH] and [Wa1] all appeal to the
positive energy theorem for asymptotically hyperbolic spaces (see
[CH], [Wa3] and some early references therein). Our proof of the
uniqueness of the AdS space-time instead appeals to the classic
positive mass theorem for asymptotically flat spaces. By this
approach we may use the classic positive mass theorem of Schoen
and Yau [SY] to drop the spin structure assumption in dimension
less than $7$.

The anti-de Sitter space-time in $n+1$ dimensions is given by
$(R^{n+1}, g_{AdS})$ where
$$
g_{AdS} = - (1+r^2)dt^2 + \frac 1{1+r^2} dr^2 + r^2 d\sigma_0 \tag
1.1
$$
in coordinates $(t, r, \theta)\in R\times [0, \infty)\times
S^{n-1}$ and $d\sigma_0$ is the standard round metric on a unit
$(n-1)$-sphere. It is a static solution to the vacuum Einstein
equation
$$
\text{Ric} - \frac 12 R \, g + \Lambda g = 0 \tag 1.2
$$
with negative cosmological constant $\Lambda = -\frac 12 n(n-1)$.
The staticity means that $(R^{n+1}, g_{AdS})$ can be constructed
by a triple $(R^n, g_H, \sqrt {1+r^2})$ where $g_H$ is the
hyperbolic metric on $R^n$ and
$$
\nabla^2 \sqrt {1+r^2} = \sqrt {1+r^2} \ g_H \tag 1.3
$$
on the hyperbolic space $(R^n, g_H)$. The simplest examples of
space-times that are asymptotically the same as the anti-de Sitter
space-time at the infinity are the so-called Schwarzschild-AdS
space-times whose metrics are given by
$$
g^+_M = - (1+r^2 - \frac {M}{r^{n-2}}) dt^2 + \frac 1{1+r^2 -
\frac {M}{r^{n-2}}} dr^2 + r^2 d\sigma_0. \tag 1.4
$$
They also satisfy the vacuum Einstein equation (1.2), but the
difference is that on the AdS space-time there is an everywhere
time-like Killing field $\frac{\partial}{\partial t}$ while this
is not so on the Schwarzschild-AdS space-times. In other words,
the AdS space-time is strictly stationary, but the
Schwarzschild-AdS space-times are not.

We will follow the idea in [AM] to give a definition of
asymptotically AdS space-times (see Definition 2.1). One can find
a good discussion of the comparisons of different definitions of
asymptotically AdS space-times in [CS]. Then we show

\proclaim{Theorem 1.1} Suppose that $(Y^{n+1}, g)$ is a strictly
stationary asymptotically AdS space-time. And suppose that $g$
satisfies the vacuum Einstein equation with negative cosmological
constant. Then $(Y^{n+1}, g)$ is static, i.e.
$$
\left\{ \aligned Y^{n+1} & = R\times \Sigma \\
g & = -Vdt^2 + h\endaligned\right. \tag 1.5
$$
where $V >0$ on $\Sigma$ and
$$
\left\{\aligned \Delta \sqrt V  & = n \sqrt V \quad  \\
Ric[h] + nh & = (\sqrt V)^{-1}\nabla^2 \sqrt V
\endaligned\right.
\tag 1.6
$$
on the Riemannian manifold $(\Sigma, h)$.
\endproclaim

By the Frobenius Theorem, staticity  is locally equivalent to
$\theta = \omega\wedge d\omega = 0$ where $\omega$ is the dual of
the given Killing vector field $X$. Instead of using topological
assumptions to write $*\theta = d \psi$ to prove the vanishing of
$\theta$ in classic Lichnerowicz argument (see [BGH], [Ca]) we
observe that
$$
d(\frac 1V \omega) = - \frac 1{V^2} i_X \theta, \tag 1.7
$$
which allows us to calculate the boundary integral to show the
vanishing of $\theta$ from the behavior of $X$ near the infinity.
We adopt the method of Fefferman and Graham [FG], [G] to construct
a preferable coordinate system near the infinity which allows us
to know the asymptotic behavior of both the metric $g$ and the
Killing field $X$ near the infinity rather precisely.

Our next goal is to prove the static solution $(\Sigma, h, \sqrt
V)$ must be the same as $(R^n, g_H, \sqrt{1+r^2})$ for some choice
of coordinates. Namely,

\proclaim{Theorem 1.2} Suppose that an asymptotically AdS
space-time is a static solution satisfying (1.5) and (1.6). Then
$(\Sigma, h, \sqrt V) = (R^n, g_H, \sqrt{1+r^2})$ for some choice
of coordinates in dimension between $3$ and $7$.
\endproclaim

Our approach is similar to the one used to prove the uniqueness of
conformally compact Einstein manifolds in [Q]. We use the global
defining function $(\sqrt V +1)^{-1}$ to turn $(\Sigma, h)$ into a
compact manifold $(\Sigma, \bar h)$ which has the round sphere as
its totally umbilical boundary and whose scalar curvature is
nonnegative. The nonnegativity of the scalar curvature follows
from the application of the strong maximum principle and the
following Bochner formula:

\proclaim{Lemma 1.3} Suppose that an asymptotically AdS space-time
is a static solution satisfying (1.5) and (1.6). Then
$$
-\Delta (V - |\nabla\sqrt V|^2 - 1) = 2 |\nabla^2\sqrt V - \sqrt V
h|^2 - \frac {\nabla\sqrt V}{\sqrt V} \cdot \nabla(V -
|\nabla\sqrt V|^2 - 1). \tag 1.8
$$
\endproclaim

Then we appeal to the recent work in [Mi](see also [ST]) to
conclude that it has to be scalar flat, which implies
$$
\nabla^2\sqrt V = \sqrt V h. \tag 1.9
$$
Note that [Mi] relies on the classic positive mass theorem of
Scoen and Yau [SY].  Theorem 1.2 then follows from the following
lemma similar to a theorem of Obata in [Ob].

\proclaim{Lemma 1.4} Suppose that $(M, g)$ is a complete
Riemannian manifold. And suppose that there is a positive function
$\phi$ such that $\nabla^2\phi = \phi g$. Then $(M, g)$ is
isometric to $(R^n, g_H)$.
\endproclaim

\head 2. Asymptotically AdS space-times \endhead

In this section we will start with an intrinsic definition of
asymptotically AdS space-times and derive some properties of a
strictly stationary asymptotically AdS space-time. We will then
prove a lemma of Lichnerowicz type similar to the one in [BGH]. We
note that, in fact, it was asked whether the uniqueness theorem in
their paper [BGH] still holds if one uses the definition of
asymptotically AdS space-times proposed by Ashtekar and Magnon in
[AM] (see also [Ha]).

Let us first introduce the AdS space-time in general dimensions.
The anti-de Sitter space-time in $(n+1)$ dimensions can be given
by $(R^{n+1}, g_{AdS})$ where
$$
g_{AdS} = - (1+r^2)dt^2 + \frac 1{1+r^2}dr^2 + r^2 d\sigma_0 \tag
2.1
$$
$d\sigma_0$ is the unit round metric on $S^{n-1}$, $t\in (-\infty,
+\infty)$, and $r\in [0, +\infty)$. In the following we will adopt
the definition of an asymptotically AdS space-time given by
Ashtekar and Magnon in [AM]. To do so, let us first discuss what a
conformal completion for a space-time is following the idea of
Penrose in [Pe]. Suppose that $Y^{n+1}$ is a manifold with
boundary $\partial Y^{n+1} = X^n$. Then $\Omega$ is said to be a
defining function of $X^n$ in $Y^{n+1}$ if

a) $\Omega > 0$, in $Y^{n+1}$;

b) $\Omega = 0$, on $X^n$; and

c) $d\Omega \neq 0$ on $X^n$. \newline A space-time $(Y^{n+1}, g)$
has a $C^k$ conformal completion if $Y^{n+1}$ is a manifold with
boundary $X^n$ and the metric $\Omega^2 g$ for a defining function
$\Omega$ of $X^n$ in $Y^{n+1}$ extends in $C^k$ to the closure of
$Y^{n+1}$.

\proclaim{Definition 2.1} A space-time $(Y^{n+1}, g)$ of dimension
$(n+1)$ is said to be asymptotically AdS if

1) $(Y^{n+1}, g)$ has a $C^k$ conformal completion and its
boundary $\partial Y^{n+1} = X^n$ is topologically $R\times
S^{n-1}$;

2) the space-time $(Y^{n+1}, g)$ satisfies the Einstein equation
with a negative cosmological constant $\Lambda$
$$
R_{ab} - \frac 12 Rg_{ab} + \Lambda g_{ab} = 8\pi G T_{ab} \tag
2.2
$$
where $\Omega^{-n} T_{ab}$ admits a $C^k$ extension to the closure
of $Y^{n+1}$;

3) $(X^n, \Omega^2 g|_{TX^n})$ is conformal to $(R\times S^{n-1},
g_0)$ where $g_0 = -dt^2 + d\sigma_0$.
\endproclaim

For the convenience, from now on, we will always assume
$$
\Lambda = -\frac {n(n-1)}2. \tag 2.3
$$

First, by definition, defining functions for $X^n$ in $Y^{n+1}$
are not unique and two different defining functions differ by a
positive function on the closure of $Y^{n+1}$. Therefore only the
class of quadratic forms $\Omega^2 g|_{TX^n}$ up to a conformal
factor is determined by $g$. Second, by requiring the fall-off of
the energy-stress tensor, one can compute that the sectional
curvature of $g$ would asymptotically go to $
-|d\Omega|^2_{\Omega^2 g}$, and conclude
$$
|d\Omega|^2_{\Omega^2 g} |_{X^n} = 1 > 0.
$$
Therefore $X^n$ is a time-like hypersurface in $(Y^{n+1}, \Omega^2
g)$. Finally the conformal flatness of the Lorentz metric
$\Omega^2 g|_{TX^n}$ depends only on the Lorentz metric $g$. In
fact, Hawking in [Ha] had already observed that locally conformal
flatness is an appropriate boundary condition.

We next want to choose a coordinate system near the boundary for
an asymptotically AdS space-time. What we will do is mostly an
analogue to the Euclidean cases which have been established in
[FG], [G]. First, we construct a special defining function, at
least in a tubular neighborhood of the boundary for each given
metric in the class $[-dt^2 + d\sigma_0]$ on the boundary by
solving a first order PDE. Namely,

\proclaim{Lemma 2.1} Suppose $(Y^{n+1}, g)$ is an asymptotically
AdS space-time, and $\Omega$ is a defining function. Then, for
each metric $\hat g = e^{2\phi} g_0$ where $g_0 = -dt^2 +
d\sigma_0$, there is a unique defining function $s$ in a tubular
neighborhood of the boundary $X^n$ in $Y^{n+1}$ such that

a) $s^2g|_{TX^n} = \hat g$;

b) $|ds|_{s^2 g} = 1$ in the tubular neighborhood.
\endproclaim
\demo{Proof} Set $s = e^w\Omega$. Then
$$
ds = e^w(d\Omega + \Omega dw)
$$
and
$$
\aligned |ds|^2_{s^2 g} & = |ds|^2_{e^{2w} \Omega^2 g} \\
& = e^{-2w}|ds|^2_{\Omega^2 g} \\
& = |d\Omega + \Omega dw|^2_{\Omega^2 g} \\
& = |d\Omega|^2_{\Omega^2 g} + 2 \Omega (d\Omega, dw)_{\Omega^2 g}
+ \Omega^2 |dw|^2_{\Omega^2 g}.
\endaligned
$$
Thus, the requirement $|ds|_{s^2 g} = 1$ is equivalent to solving
$$
2(d\Omega, dw)_{\Omega^2 g} + \Omega |dw|^2_{\Omega^2 g} = \frac
{1 -|d\Omega|^2_{\Omega^2 g}}\Omega. \tag 2.4
$$
The boundary condition is determined as follows: if we denote
$\Omega^2 g |_{TX^n} = e^{2\psi}g_0$, then
$$
w|_{X^n}  = \phi - \psi. \tag 2.5
$$
It is easily seen that (2.4) and (2.5) is non-characteristic.
Notice that both $d\Omega$ and $dw$ are space-like.
\enddemo

\proclaim{Lemma 2.2} Suppose $(Y^{n+1}, g)$ is an asymptotically
AdS space-time. Suppose that $s$ is the special defining function
obtained in Lemma 2.1 for which $s^2 g|_{TX^n} = g_0$. Then
$$
g = s^{-2}(d s^2 + g_s) \tag 2.6
$$
where
$$
g_s = -(1+\frac {s^2}4)^2dt^2 + (1 - \frac {s^2}4)^2 d\sigma_0  +
O(s^n). \tag 2.7
$$
\endproclaim
\demo{Proof} The proof again is adopted from the argument given in
[FG], [G]. By the fall-off condition of the energy-momentum tensor
$T_{ab}$ one can rewrite the equation (2.2) in coordinates
$R\times S^{n-1}\times [0, \epsilon)$ near the boundary as
$$
h_{ab}'' + (1-n) h_{ab}' - h^{cd}h_{cd}'h_{cd} - s
h^{cd}h_{ac}'h_{bd}' + \frac 12 s h^{cd}h_{cd}' h_{ab}' -2s
R_{ab}[h] = O(s^n). \tag 2.8
$$
where $h$ stands for $g_s$ for convenience. The signature of $g_s$
here does not make any difference in terms of solving the
expansion of $g_s$. Therefore, similar to what is known for
Euclidean case, all odd order terms of order $\leq n-1$ vanish and
all even order terms of order $\leq n-1$ is determined by the
metric $g_0$ on $X^n$. Moreover, when $n$ is odd, the $n$th order
term is traceless; when $n$ is even, in general one would need to
add one more term in the order of $s^n\log s$ which is traceless
and determined by $g_0$ while the trace part of the $n$th order is
also determined by $g_0$. By comparing to the AdS space
$$
g_{AdS} = s^{-2}(d s^2 - (1+\frac {s^2}4)^2dt^2 + (1-\frac
{s^2}4)^2d\sigma_0)
$$
which is of the same boundary metric $-dt^2 + d\sigma_0$, we may
complete the proof.
\enddemo

\proclaim{Remark 2.3} In the above argument, it is clear that a
weaker fall-off condition of the energy-momentum would imply a
weaker control of the asymptotic of the metric $g$.
\endproclaim

Next we will follow [BGH] to restrict ourselves to the so-called
strictly stationary space-time. That is to assume, for an
asymptotically AdS space-time, there is a global everywhere
time-like Killing field which approaches $\frac{\partial}{\partial
t}$ asymptotically towards the boundary. In [BGH] it was shown
that a strictly stationary asymptotically AdS space-time (by their
definition in dimension 3) which solves the vacuum Einstein
equations with negative cosmological constant must be a static
one. Before proceeding to prove the staticity we want to study the
asymptotic behavior of a Killing field that approaches
$\frac{\partial}{\partial t}$ at the infinity. We will use the
favorable coordinates constructed in Lemma 2.2. Denote the Killing
field by
$$
X = a (s, t, \sigma)\frac{\partial}{\partial s} + b(s,t,\sigma)
\frac{\partial}{\partial t} +
c^i(s,t,\sigma)\frac{\partial}{\partial \theta_i}, \tag 2.9
$$
where $(s, t, \theta_1, \cdots, \theta_{n-1})\in [0,
\epsilon)\times R\times S^{n-1}$. For similar computations, please
see [Wa2]. First of all, by the boundary condition, we know that
$$
b(0, t, \theta) = 1, \quad c^i(0, t, \theta) = 0, \forall i= 1,
\cdots, n-1. \tag 2.10
$$
Computing $X g(\frac{\partial}{\partial s},
\frac{\partial}{\partial s})$ one gets
$$
\frac{\partial a}{\partial s} = \frac as. \tag 2.11
$$
Therefore $a(s,t,\theta) = sa(t, \theta)$. For convenience we
denote by $ t=\theta_0$ and $b = c^0$,  and use Greek letters to
include zero. Computing $Xg( \frac{\partial}{\partial s},
\frac{\partial}{\partial \theta_\alpha})$ one gets
$$
g_s(\frac{\partial}{\partial \theta_\alpha},
\frac{\partial}{\partial \theta_\beta})\frac {\partial
c^\beta}{\partial s} +  s \frac {\partial a}{\partial
\theta_\alpha} = 0.
$$
Therefore
$$
\left\{\aligned b(s,t,\theta) & = 1 - \int_0^s ug_u^{0\beta}du \,
\frac {\partial a}{\partial \theta_\beta} \\
c^\alpha (s, t, \theta) & = - \int_0^s ug_u^{\alpha\beta}du \,
\frac {\partial a}{\partial \theta_\beta}. \endaligned\right.\tag
2.12
$$
In the other directions one gets
$$
g_s(\frac {\partial}{\partial \theta_\alpha}, \frac
{\partial}{\partial \theta_\gamma})c^\gamma_{\ , \beta} +
g_s(\frac {\partial}{\partial \theta_\beta}, \frac
{\partial}{\partial \theta_\gamma})c^\gamma_{\ , \alpha} = 0. \tag
2.13
$$
Notice that, if we denote the Christoffel symbols of metric $g_s$
on the slices by $\Gamma^\alpha_{\beta\gamma}$ and $\bar
\Gamma^a_{bc}$ for the ones of metric $g$, then
$$
\Gamma^\alpha_{\beta\gamma} = \bar \Gamma^\alpha_{\beta\gamma},
\quad \bar \Gamma^\alpha_{s\beta} = \frac 12 h^{\alpha\gamma}\frac
{\partial}{\partial s} h_{\gamma\beta} - \frac 1s
\delta_{\alpha\beta}, \tag 2.14
$$
where again, for convenience, we use $h=g_s$. Thus
$$
h_{\alpha\gamma}(\frac{\partial c^\gamma}{\partial \theta_\beta} -
c^\delta\Gamma^\gamma_{\delta\beta}) +
h_{\beta\gamma}(\frac{\partial c^\gamma }{\partial \theta_\alpha}-
c^\delta\Gamma^\gamma_{\delta\alpha}) = 2 a h_{\alpha\beta} - sa
\frac {\partial}{\partial s}h_{\alpha\beta}. \tag 2.15
$$
Taking $s=0$ in (2.15) we immediately see that $a(t,\theta)=0$ and
surprisingly get $X = \frac{\partial}{\partial t}$ in this
neighborhood. Moreover (2.15) implies $g_s$ is independent of $t$
in this neighborhood. Let us summarize what we obtained in a
lemma.

\proclaim{Lemma 2.4} Suppose that $(Y^{n+1}, g)$ is an
asymptotically AdS space-time and that $s$ is the special defining
function such that $s^2 g|_{TX} = -dt^2 + d\sigma_0$. If $X$ is a
Killing field that approaches to $\frac {\partial}{\partial t}$ at
the infinity, then $X = \frac {\partial}{\partial t}$ and $g_s$ is
independent on $t$ in a tubular neighborhood of the boundary.
\endproclaim

A strictly stationary asymptotically AdS space-time $(Y^{n+1}, g)$
comes with a proper and free $R^1$ action. Here properness of the
action comes from the causality axiom (cf. [Ca]). Therefore, by
Theorem 1.11.4 in [DK], we know that $Y$ is a principle bundle
over the smooth orbit space $Y/R$. Thus topologically, $Y = R
\times M$ for some smooth $n$-manifold $M$ with boundary $S^{n-1}$
in the light of the discussion above on the asymptotic behavior of
the Killing field $X$.

Now let us discuss the staticity of a space-time. A good reference
for this discussion is [Ca]. A space-time is said to be static if
there is an everywhere time-like Killing field whose trajectories
are everywhere orthogonal to a family of space-like hypersurfaces.
Let us introduce some notations. Let $\{e_a\}$ be an orthonormal
frame and $\{w^a\}$ be its co-frame. Suppose $X = k^ae_a$ is an
everywhere time-like Killing field and let $\omega = k_aw^a$. In
this notation $X$ is a Killing field if
$$
k_{a, b} + k_{b, a} = 0. \tag 2.16
$$
By Frobenius Theorem, staticity is equivalent to asking that the
differential ideal generated by the differential $\omega$ be
closed under exterior differentiation, i.e.
$$
\theta = \omega \wedge d\omega = 0. \tag 2.17
$$
A connected static space-time becomes $R\times \Sigma$ where
$\Sigma$ is a static slice (topologically the same as $M$), and
the metric $g = -Vdt^2 + g_\Sigma$ where $V = - k^ak_a$ and
$g_\Sigma$ is metric on the slice $\Sigma$ of Euclidean signature.
The following lemma of Lichnerowicz type is a straightforward
generalization of a lemma in [BGH] (see also [Ca] for more
details), but the proof is adapted for general dimension and uses
no additional topological assumptions. Notice that the definition
of an asymptotically AdS space-time in [BGH] is different from
ours in this note.

\proclaim{Lemma 2.5} Any strictly stationary asymptotically AdS
space-time  $(Y^{n+1}, g)$ \newline which satisfies the vacuum
Einstein equations with negative cosmological constant $\Lambda$
is static.
\endproclaim
To prove this lemma we observe

\proclaim{Lemma 2.6}
$$
\frac {i_X \theta}{V^2} = - d(\frac \omega V) \tag 2.18
$$
where
$$
V = - k^ak_a.
$$
\endproclaim
\demo{Proof} We simply compute
$$
d(\frac \omega V) = \frac {d\omega}V - \frac 1{V^2}dV\wedge \omega
= -\frac 1{V^2}(d\omega (-V) + dV\wedge \omega) = -\frac {i_X
\theta}{V^2}.
$$
Because
$$
i_X \omega = \omega(X) = k_ak^a = - V
$$
and
$$
\aligned i_X dw & = i_X (k_{a,b}w^b\wedge w^a) = k_{a,b}i_X w^b
w^a- k_{a,b}w^b i_X w^a \\
& = k_{a,b}k^bw^a - k_{a, b}k^a w^b = -2k^bk_{b,a}w^a = dV.
\endaligned
$$
\enddemo

\demo{Proof of Lemma 2.5} Let us consider the Hodge dual $*\theta$
of $\theta$ . Since
$$
d^*\theta = (k_{[a,b}k_{c]})^c w^a\wedge w^b = \frac 23 k^c
R_{c[a}k_{b]} w^a\wedge w^b = 0 \tag 2.19
$$
due to the fact that $R_{ab} = n \eta_{ab}$ where $\eta_{ab}$ is
the standard Minkowski metric (please see Chapter 6 in Part II of
[Ca]) and
$$
i_X (*\theta) = *(\theta\wedge \omega) = 0, \tag 2.20
$$
it follows
$$
d (\frac \omega V \wedge *\theta)  = -\frac {i_X \theta \wedge
*\theta}{V^2} = - \frac {i_X (\theta\wedge *\theta)}{V^2}. \tag
2.21
$$
The next step is to integrate over a space-like hypersurface
$\Sigma_\epsilon$ whose boundary is a large $(n-1)$-sphere $
S^{n-1}_\epsilon = \{ s=\epsilon, t = c\}$ in the preferable
coordinates. We therefore have
$$
\int_{\Sigma_\epsilon} \frac {|\theta|^2 \omega(N)}{V^2} d\sigma =
- \int_{S^{n-1}_\epsilon} \frac \omega V\wedge *\theta \tag 2.22
$$
where $N$ is the unit normal of $\Sigma_\epsilon$, and $d\sigma$
is the volume element of $\Sigma_\epsilon$ in the space-time.
Notice that $\omega(N)> 0$ and $\theta$ is space-like since
$\theta\wedge\omega=0$. Now let us recall that the Killing field
$X$ is just $\frac{\partial}{\partial t}$ in the preferable
coordinates. Thus
$$
\omega |_{S^{n-1}_\epsilon}= - V dt + g_{0k}d\theta^k = -V dt +
s^{n-2}\tau_{0k}d\theta^k = \epsilon^{n-2} \tau_{0k}d\theta^k \tag
2.23
$$
where
$$
g = s^{-2}(ds^2 + (1+\frac {s^2}4)^2dt^2 + (1-\frac
{s^2}4)^2d\sigma_0 + s^n\tau) \tag 2.24
$$
and
$$
V = s^{-2}(1+\frac {s^2}4)^2 - s^{n-2}\tau_{00}. \tag 2.25
$$
Then
$$
d\omega = -dV\wedge dt + d g_{0k}\wedge d\theta^k,
$$
$$
\theta = - Vdt\wedge dg_{0k}\wedge d\theta^k -
s^{n-2}\tau_{0k}d\theta^k\wedge dV\wedge dt + s^{n-2}\tau_{0k}
d\theta^k\wedge dg_{0l}\wedge \theta^l, \tag 2.26
$$
and
$$
\aligned *\theta|_{S^{n-1}_\epsilon} & =  (-V \frac {\partial
g_{0k}}{\partial s}  + g_{0k}\frac {\partial V}{\partial
s})*(dt\wedge ds\wedge d\theta^k) \\
& =  C s^{n-5} *(dt\wedge ds\wedge d\theta^k) \\
& = C d\theta^1\wedge\cdots \widehat{d\theta^k}\wedge\cdots\wedge
d\theta^{n-1}
\endaligned \tag 2.27
$$
where $C$ stands for some function on $S^{n-1}$. Therefore
$$
\frac \omega V \wedge *\theta |_{S^{n-1}_\epsilon} = O(\epsilon^n)
d\theta^1\wedge d\theta^2 \wedge\cdots\wedge d\theta^{n-1}. \tag
2.28
$$
This implies that $\theta = 0$ on the hypersurface $\Sigma$. But
$\Sigma$ is arbitrary, so $\theta = 0$ on $Y^{n+1}$, which
finishes the proof.
\enddemo

Let us conclude this section by making it clear what a static
asymptotically AdS space-time which satisfies the vacuum Einstein
equations with negative cosmological constant $\Lambda$ is. We
state our observation in the following lemma.

\proclaim{Lemma 2.7} Under the assumption of Lemma 2.5, in the
preferable coordinate system at the infinity, indeed, a slice of
constant $t$ is a static slice, i.e. $\frac{\partial}{\partial t}$
is orthogonal to the slice of constant $t$.
\endproclaim
\demo{Proof} Consider the conformal completion $(Y^{n+1}, \bar g)$
where $\bar g = ds^2 + g_s$. By the construction of the preferable
coordinate system, each curve $\gamma(s) = (s, t, \theta_1,
\cdots, \theta_{n-1})$ is a geodesic from the point $(0, t,
\theta_1, \cdots, \theta_{n-1})$ in the space-time $(Y^{n+1}, \bar
g)$. On the other hand, a static slice $\Sigma$ of $(Y^{n+1}, g)$
is still a maximum integral hypersurface which is orthogonal to
$\frac{\partial}{\partial t}$ everywhere with respect to $\bar g$.
Because $\bar g = ds^2 + g_s$ is independent of $t$, such $\Sigma$
is totally geodesic in $(Y^{n+1}, \bar g)$. Therefore a geodesic
emanating from a boundary point $(0, t_0, \theta_1, \cdots,
\theta_{n-1})$ with respect to the metric $\bar g$ stays in a
static slice. Thus a slice of constant $t$ coincides with a static
slice. So the proof is complete.
\enddemo

We summarize our result in the following theorem:

\proclaim{Theorem 2.8} Suppose $(Y^{n+1}, g)$ is a strictly
stationary asymptotically AdS space-time that satisfies the vacuum
Einstein equations with negative cosmological constant $\Lambda =
-\frac {n(n-1)}2$. Then $Y^{n+1}=R\times \Sigma^n$,
$$
g= - Vdt^2 + h \tag 2.29
$$
and, on $\Sigma$,
$$
\left\{\aligned \Delta \sqrt V  & = n \sqrt V \quad  \\
Ric[h] + nh & = (\sqrt V)^{-1}\nabla^2 \sqrt V ,
\endaligned\right.
\tag 2.30
$$
where $h$ is the metric induced from $g$ on a static slice
$\Sigma$ of Euclidean signature. Moreover $(\Sigma, h)$ is
conformally compact of the same regularity as of the conformal
completion of $(Y^{n+1}, g)$ with the conformal infinity
$(S^{n-1}, [d\sigma_0])$ where
$$
V^{-1} h|_{TS^{n-1}} = d\sigma_0.
$$
\endproclaim

\head 3. Static asymptotically AdS space-times \endhead

In this section we study static asymptotically AdS space-times. We
will prove the uniqueness of static asymptotically AdS
space-times. In dimension $3+1$, with a bit restrictive definition
of asymptotically AdS space-times, the uniqueness was first proved
in [BGH] (see also [CS]). Then, assuming spin structure for $n >
3$, the uniqueness of static solutions $(M^n, g, V)$ to the vacuum
Einstein equations with negative cosmological constant was
established in [Wa1](see the definition of a static solution $(M,
g, V)$ in [CS], [Wa1]). Our proof will not use the spin structure
in dimensions higher than three, but instead will rely on a recent
work of Miao [Mi] (see also [ST]) which in turn depends on the
classic positive theorem of Schoen and Yau [SY] for asymptotically
flat manifolds.

By Theorem 2.8 in the previous section, a static asymptotically
AdS space-time which satisfies the vacuum Einstein equations with
negative cosmological constant is given by a static solution
$(\Sigma, h, \sqrt V)$ in our notation. Therefore, by Lemma 2.2 in
the previous section, we know that
$$
h = s^{-2}(ds^2 +(1 - \frac {s^2}4)^2 d\sigma_0 + \tau s^n +
o(s^n)), \tag 3.1
$$
$$
V = s^{-2}((1 + \frac {s^2}4)^2 - \alpha s^n + o(s^n)), \tag 3.2
$$
and
$$
V - |\nabla\sqrt V|^2 - 1 = n\alpha s^{n-2} + o(s^{n-2}) \tag 3.3
$$
where $\alpha = - \text{Tr}_{d\sigma_0}\tau$ (these were known in
[Wa1]).

To motivate our argument in this section we recall the following
fact about the static solution $(B^n, (\frac 2{1-|x|^2})^2 |dx|^2,
\frac {1+|x|^2}{1-|x|^2})$ associated with the AdS space-time
$(R^{n+1}, g_{AdS})$. Namely, if one uses the global defining
function
$$
u = \frac 1{\sqrt V + 1}= \frac {1-|x|^2}2,
$$
then
$$
u^2 h = |dx|^2.
$$
Therefore, for a static solution $(\Sigma, h, \sqrt V)$, if we
denote $u = \frac 1{\sqrt V + 1}$, then $u$ is a global defining
function for $S^{n-1}$ in $\Sigma$ and
$$
u^2 h = \frac 1{s^2(\sqrt V+1)^2}(ds^2 + (1-\frac
{s^2}4)^2d\sigma_0 + \tau s^n + o(s^n))
$$
where
$$
s^2(\sqrt V+1)^2 = (\sqrt {s^2V} +s)^2 = s^2V + 2s\sqrt {s^2V}
+s^2 = 1 + 2s + O(s^2).
$$
So
$$
u^2h = (1 + 2s + O(s^2))ds^2 + (1 + 2s + O(s^2))d\sigma_0 +
O(s^2). \tag 3.4
$$
Thus $(\Sigma, u^2h)$ is a compact manifold with the standard
$(n-1)$-sphere as its boundary and the second fundamental form for
$\partial\Sigma$ in $\Sigma$ is $d\sigma_0$ (i.e. the boundary is
totally umbilical). In the light of (2.30) one may compute the
scalar curvature for $u^2h$ as follows:
$$
\aligned R[u^2h] & = u^{-\frac {n+2}2}( - \frac
{4(n-1)}{n-2}\Delta u^{\frac {n-2}2} + R[h]u^{\frac {n-2}2})\\
& = n(n-1) ( V - |\nabla\sqrt V|^2 - 1),
\endaligned\tag 3.5
$$
which goes to zero as $s\to 0$ by (3.3). We observe the following
lemma which will allow us to apply the Strong maximum principle to
conclude that $R[u^2h] \geq 0$. Namely,

\proclaim{Lemma 3.1}
$$
-\Delta (V - |\nabla\sqrt V|^2 - 1) = 2 |\nabla^2\sqrt V - \sqrt V
h|^2 - \frac {\nabla\sqrt V}{\sqrt V} \cdot \nabla(V -
|\nabla\sqrt V|^2 - 1). \tag 3.6
$$
\endproclaim
\demo{Proof} We simply compute
$$
-\Delta V = 2(-nV - |\nabla\sqrt V|^2)
$$
and
$$
\aligned \Delta |\nabla\sqrt V|^2 & = 2 (|\nabla^2 \sqrt V|^2 +
(\sqrt V)_i (\sqrt V)_{ijj}) \\
& = 2 (|\nabla^2 \sqrt V|^2 +  (\sqrt V)^{-1}(\sqrt V)_i(\sqrt
V)_{ij}(\sqrt V)_j)).
\endaligned
$$
Therefore
$$
\aligned -\Delta & (V  - |\nabla\sqrt V|^2 - 1)  \\
& = 2 |\nabla^2\sqrt V - \sqrt V h|^2 - \frac {(\sqrt V)_i}{\sqrt
V}(2(\sqrt V)(\sqrt V)_i - 2(\sqrt V)_{ij}(\sqrt
V)_j) \\
& = 2 |\nabla^2\sqrt V - \sqrt V h|^2 - \frac {\nabla\sqrt
V}{\sqrt V} \cdot \nabla(V - |\nabla\sqrt V|^2 - 1).
\endaligned
\tag 3.7
$$
\enddemo

\proclaim{Theorem 3.2} Suppose that $(\Sigma, h, \sqrt V)$ is a
static solution to the vacuum Einstein equations with negative
cosmological constant, i.e. $(\Sigma, h, \sqrt V)$ comes from
Theorem 2.8 in the previous section. Then $(\Sigma, h, \sqrt V)=
(B^n, (\frac 2{1-|x|^2})^2|dx|^2, \frac {1+|x|^2}{1-|x|^2})$ for
some choice of coordinate in dimension between $3$ and $7$.
\endproclaim
\demo{Proof} We consider the defining function $u = \frac 1{\sqrt
V + 1}$ and the compact manifold $(\Sigma, u^2 h)$. From (3.4) we
know that $\partial\Sigma = S^{n-1}$ and $u^2h|_{\partial\Sigma} =
d\sigma_0$. Moreover, from (3.4), we know that the standard round
$S^{n-1}$ is the boundary of $(\Sigma, u^2h)$ and has the second
fundamental form $d\sigma_0$ (i.e. it is totally umbilical). On
the other hand, by (3.5) and (3.3), the scalar curvature $R[u^2h]$
goes to zero as $s\to 0$. Using the above Lemma 3.1 and the strong
maximum principle, we therefore conclude that $R[u^2h] \geq 0$.

Now we appeal to the recent work of Miao [Mi] (see also works of
Shi and Tam [ST]). We apply the work in [Mi] to the manifold $(M,
\Cal G)$ where $\Omega = \Sigma$ and $g_- = u^2h$, and $M\setminus
\Omega = R^n\setminus B^n$ and $g_+$ is the Euclidean metric. By
Corollary 5.1 in [Mi], for example, we then conclude that $R[u^2h]
\equiv 0$ for $ 3\leq n\leq 7$.

Back to $(\Sigma, h, \sqrt V)$, in the light of Lemma 3.1 again,
we observe that
$$
\nabla^2 \sqrt V = \sqrt V h. \tag 3.8
$$
Similar to what was proved in [Ob], we prove that (3.8) implies
that $(\Sigma, h)$ is isometric to the standard hyperbolic space
form and $V = \sqrt {1+r^2}$ for some choice of coordinates in the
following lemma. Then the proof of theorem is complete.
\enddemo

\proclaim{Lemma 3.3} Suppose that $(M^n, g)$ is a complete
Riemannian manifold. And suppose that there is a positive function
$\phi$ such that
$$
\nabla^2 \phi = \phi g. \tag 3.9
$$
Then $(M^n, g)=(R^n, g_H)$ and $\phi = c \sqrt {1+r^2}$ for some
choice of coordinates.
\endproclaim
\demo{Proof} First one observes that, $\phi$ has one and only one
global minimum point $p_0$ on $M$. Due to the homogeneity of
(3.9), one may assume that $\phi (p_0) = 1$. Let us consider a
geodesic $\gamma (s)$ emanating from $p_0$ and parameterized with
its length $s$. Then, along this geodesic, for $\phi (\gamma(s))$,
we have
$$
\left\{ \aligned \phi'' - \phi & = 0\\
\phi (0) & = 1\\
\phi'(0) & = 0.\endaligned\right. \tag 3.10
$$
Thus $\phi (s) = \cosh s$. Now take an othonormal base $X^0= \frac
{\partial}{\partial s}, X^1, X^2,\cdots X^{n-1}$ at $p_0$ and
parallel translate them along $\gamma (s)$. We want to calculate
$d(\text{exp}_{p_0})_{(sX^0)}(X^k)$. That is, we compute the
Jacobi field $Y^k(s)$ along $\gamma(s)$ such that
$$
\left\{\aligned \nabla_{X^0}\nabla_{X^0} Y^k + R(Y^k, X^0)X^0 & =
0\\  Y^k(0) & = 0 \\ \nabla_{X^0}Y^k (0) & =
X^k(0).\endaligned\right. \tag 3.11
$$
Let $Y^k (s) = \sum f_i(s) X^i(s)$. Then (3.11) becomes
$$
\left\{\aligned f_i'' X^i + f_i R_{0i0j} X^j & = 0 \\
f_i & = 0 \\ f_i' & = \delta_{ik}.\endaligned\right. \tag 3.12
$$
Notice that, by (3.9) and Ricci identity,
$$
\phi_{a,bc} - \phi_{a,cb} = \phi_d R_{dabc} = \sinh s R_{0abc} =
\phi_c \delta_{ab} - \phi_{b}\delta_{ac}, \tag 3.13
$$
which gives us
$$
R_{0i0j} = \frac 1{\sinh s} (\phi_j\delta_{i0} -
\phi_0\delta_{ij}) = - \delta_{ij}. \tag 3.14
$$
Plugging into (3.12), we have
$$
\left\{\aligned f_i'' - f_i & = 0 \\
f_i(0) & = 0\\
f_i'(0) & = \delta_{ik}.\endaligned\right.\tag 3.15
$$
Thus $Y^k(s) = \sinh s X^k(s)$.  To show that $(M, g)$ is a
hyperbolic space form, we use the exponential map $exp_{p_0}$
which takes the tangent space $T_{p_0}M$ onto $M$ in the light of
completeness. Clearly this gives a nice global coordinate chart.
Next we want to calculate the metric $g$ under these coordinates.
Let us use spherical coordinates for $T_{p_0}M$, that is, $(s,
v)\in [0, \infty)\times S^{n-1}$ and $exp_{p_0}(sv) \in M$. By the
above calculations of Jacobi fields, we immediately have
$$
g = ds^2 + (\sinh s)^2 \ d\sigma_0, \tag 3.16
$$
which is the hyperbolic metric. So $(M, g)$ is a hyperbolic space
form. Finally let us point out that, if we denote $r=\sinh s$,
then $\phi = \cosh s = \sqrt {1+r^2}$.
\enddemo

\vskip 0.1in
\noindent {\bf References}:

\roster
\item"{[ACD]}" M. Anderson, P. Chrusciel and E. Delay,
Non-trivial, static, geodesically complete, vacuum space-times
with a negative cosmological constant, arXiv: hep-th/0211006.

\vskip 0.1in
\item"{[AM]}" A. Ashtekar and A. Magnon, Asymptotically anti-de
Sitter space-times, Class. Quantum Grav. 1 (1984) L39-L40.

\vskip 0.1in
\item"{[BGH]}" W. Bocher, G.W. Gibbons and G.T. Horowitz,
Uniqueness theorem for anti-de Sitter spacetime, Phys. Review D.
(3) 30 (1984), no. 12, 2447-2451.

\vskip 0.1in
\item"{[Ca]}" B. Carter, Black hole equilibrium states, Part II,
in ``Black Holes", edited by C. DeWitt and B. DeWitt (New York,
1973).

\vskip 0.1in
\item"{[CH]}" Piotr Chru\'{s}ciel and M. Herzlic, The mass of
asymptotically hyperbolic Riemannian manifolds, Preprint
math.DG/0110035.

\vskip 0.1in
\item"{[CS]}" Piotr Chru\'{s}ciel and Walter Simon, Towards the
classification of static vacuum spacetimes with negative
cosmological constant, J. Math. Phys. {\bf 42} (2001) no.4
1779-1817.

\vskip 0.1in
\item"{[DK]}" J.Duistermaat and J. Kolk, ``Lie groups",
Springer-Verlag, Berlin, New York, 2000.

\vskip 0.1in
\item"{[FG]}" C. Fefferman, and C.R. Graham, Conformal invariants, in
{\it The mathematical heritage of Elie Cartan}, Asterisque, 1985,
95-116.

\vskip 0.1in
\item"{[GSW]}" G.J. Galloway, S. Surya and E. Woolgar, On the
geometry and mass of static, asymptotically AdS space-times,and
uniqueness of the AdS soliton, arXiv: hep-th/0204081.

\vskip 0.1in
\item"{[G]}" C. R. Graham, Volume and Area renormalizations for
conformally compact Einstein metrics. The Proceedings of the 19th
Winter School "Geometry and Physics" (Srn\`{i}, 1999). Rend. Circ.
Mat. Palermo (2) Suppl. No. 63 (2000), 31--42.

\vskip 0.1in
\item"{[Ha]}" S. Hawking, The boundary conditions for gauged
supergravity, Phys. Lett. B 126 (1983), no. 3-4, 175--177.

\vskip 0.1in
\item"{[Mi]}" P. Miao, Positive mass theorem on manifolds
admitting corners along a hypersurface, ArXiv: math-ph/0212025.

\vskip 0.1in
\item"{[Ob]}" M. Obata, Certain conditions for a Riemannian
manifold to be isometric with a sphere, J. Math. Soc. Japan. 14.
(1962) no. 3, 333-340.

\vskip 0.1in
\item"{[Pe]}" R. Penrose, Asymptotic properties of fields and
space-times, Phys. Rev. Lett. 10 (1963) 66--68.

\vskip 0.1in
\item"{[Q]}" Jie Qing, On the rigidity for conformally
compact Einstein manifolds,  International Mathematics Research
Notices, 21 (2003) 1141-1153, ArXiv: math.DG/0305084.

\vskip 0.1in
\item"{[SY]}" R. Schoen, Variational theory for the total scalar curvature
functional for Riemannian metrics and related topics. Topics in
calculus of variations (Montecatini Terme, 1987), 120--154,
Lecture Notes in Math., 1365, \newline Springer, Berlin, 1989.

\vskip 0.1in
\item"{[ST]}" Y. Shi and L. Tam, Positive mass theorem and the
boundary behaviors of a compact manifolds with nonnegative scalar
curvature, \newline arXiv: math.DG/0301047.

\vskip 0.1in
\item"{[Wa1]}" X. Wang, Uniqueness of AdS space-time in any dimension,
\newline arXiv: math.DG/0210165.

\vskip 0.1in
\item"{[Wa2]}" X. Wang, On conformally compact Einstein
manifolds, Math. Res. Lett. 8 (2001), no. 5-6, 671-688.

\vskip 0.1in
\item"{[Wa3]}" X. Wang, The mass of Asymptotically hyperbolic
manifolds, J. Diff. Geo. 57 (2001), no. 2, 273-299.

\endroster
\enddocument